\documentclass[fina2,3l,leqno]{siamltex}
\usepackage{amsmath}
\usepackage{amsfonts}
\usepackage{amssymb}
\usepackage{enumerate}
\usepackage{subfigure}
\usepackage{color}
\usepackage{graphicx,float}
\usepackage{verbatim}
\usepackage{url}
\usepackage{alltt}   
\usepackage{algorithm2e}
\usepackage{qtree}
\usepackage[english]{babel}
\usepackage{wrapfig}

\usepackage{simplemargins}
\setleftmargin{2in}
\setrightmargin{2in}
\settopmargin{2in}
\setbottommargin{2in}

\newcounter{saveeqn}%

\newtheorem{theo}{Theorem}[section]
\newtheorem{lemm}{Lemma}[section]
\newtheorem{rema}{Remark}[section]
\newtheorem{Defn}{Definition}[section]
\newtheorem{Prop}{Conjecture}[section]

\author{S.J. Gismondi\footnotemark[1]}
\title{In Honour of Ted Swart} 
\begin{document}
	


	
\maketitle
\renewcommand{\thefootnote}{\fnsymbol{footnote}}
\footnotetext[1]{\textit{University of Guelph, Canada, Email: \textup{\nocorr \texttt{gismondi@uoguelph.ca}}}}
\pagestyle{myheadings}
\markboth{Gismondi}{Ted's Tribute} 


\begin{wrapfigure}{R}{0.4\textwidth}
\centering
\includegraphics[width=38mm,height=44mm]{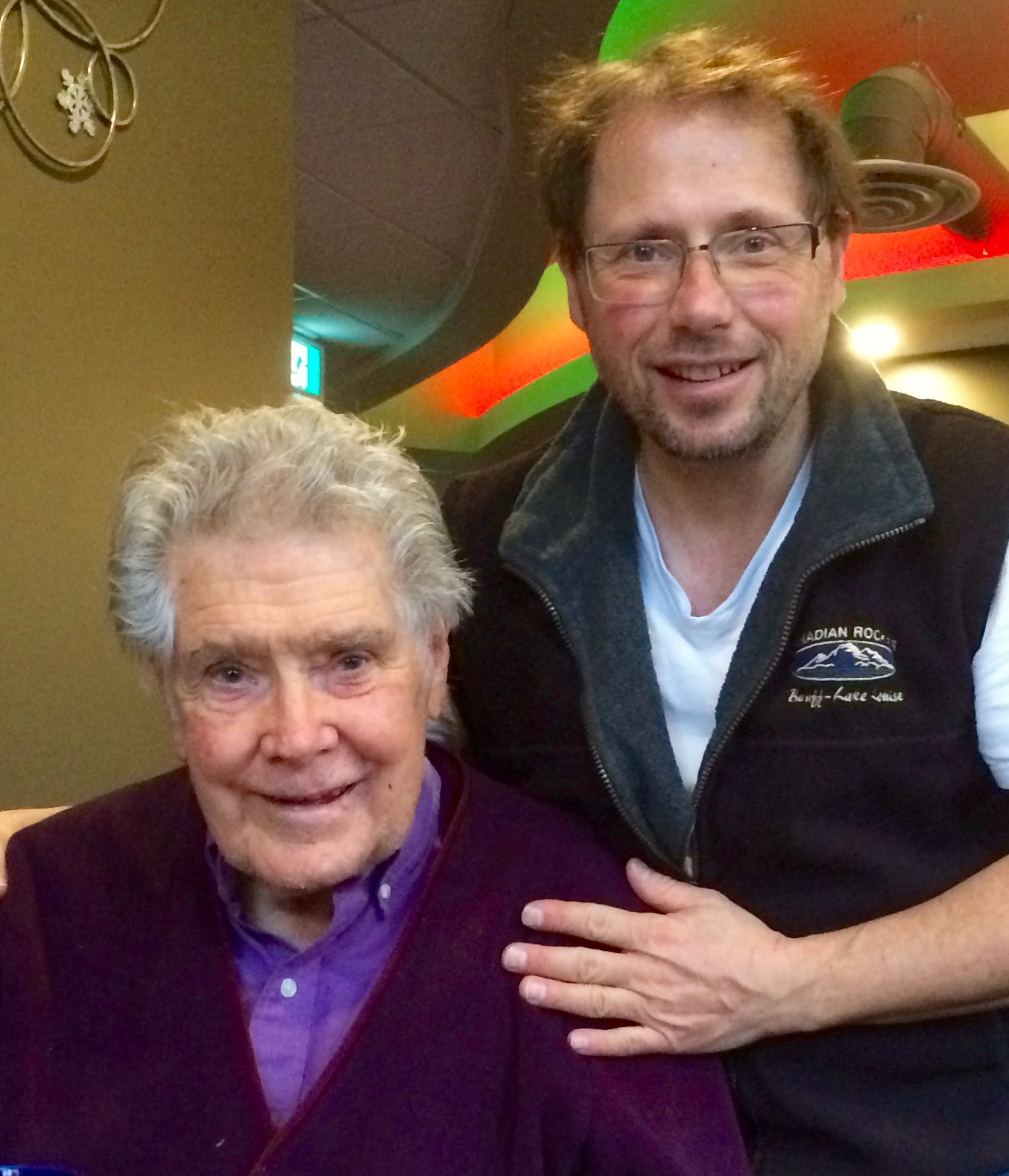}
\label{TedSteve}
\end{wrapfigure}

\footnotesize
\noindent \textbf{Forward Summary.} This is a tribute to my dear life-long friend, mentor and colleague Ted Swart\footnotemark[2]. It includes anecdotal stories and memories of our times together, and also includes a new academic contribution in his honour, Ted's polytope.

Ever since I met Ted, we loved to discuss new ideas. It didn't matter what they were. Almost always, we talked about how to model graph theory problems via LP. We'd typically meet in the afternoon, present an idea, and then challenge and encourage one another. If we were meeting at Ted's home, Ted's wife Diana would prepare snacks for us. We worked together like it was a game. We'd laugh, scribble on the board, maybe have dinner, and then go home. A few hours later, one of us would call the other with an even better idea. These ideas would survive for days, weeks and sometimes even months. But eventually we'd find a flaw and run out of energy. Ted would simply say `the truth is the truth'. But a few weeks later we'd be discussing yet another idea. This went on for 35 years, and ended only when Covid-19 restrictions made it too difficult to meet in person. Every now and then a paper might come about e.g. we came upon a polyhedral model of a \textbf{coNP-complete} problem in one of these marathon discussions. I will miss him very much.

Now about Ted's polytope, `tweeks' made to the Birkhoff polytope $B_n$ endow Ted's polytope $\mathcal{T}_n(\epsilon)$ with a special `tunable' parameter $\epsilon=\epsilon$($n$). Observe how $B_n$ can be viewed as the convex hull of both the TSP polytope, and the set of non-tour permutation extrema, and, that its extended formulation is compact. Tours (connected 2-factor `permutation matrices when viewed as adjacency matrices') can be distinguished from non-tours (disconnected 2-factor `permutation matrices') where $\epsilon$ scales the magnitude of tweeks made to $B_n$. For $\epsilon > 0$, $\mathcal{T}_n(\epsilon)$ is `tuned' so that the convex hull of extrema corresponding to transformed tours is lifted from $B_n$, and separated (by a hyperplane) from the convex hull of extrema corresponding to translated non-tours. This leads to creation of the feasible region of an LP model that can decide existence of a tour in a graph based on an extended formulation of the TSP polytope. That is, by designing for polynomial-time distinguishable tour extrema\footnotemark[3] embedded in a subspace disjoint from non-tour extrema, \textbf{NP-completeness} strongholds come into play, necessarily expressed in a non-compact extended  formulation of $\mathcal{T}_n(\epsilon)$ i.e. a compact extended formulation of the TSP polytope cannot exist. No matter, Ted would have loved these ideas, and $\mathcal{T}_n(\epsilon)$ might one day yet be useful in the study of the \textbf{P} versus \textbf{NP} conundrum\footnotemark[4]. In summary, $\mathcal{T}_n(\epsilon)$ is a perturbed $B_n$ i.e. the convex hull of both an $\epsilon$-stretched TSP polytope, and the set of translated non-tour permutation extrema i.e. a TSP-like polytope and separable non-tour extrema.

\footnotetext[2]{Ted passed in January 2023 at 94 years of age. We were the best of friends.}

\footnotetext[4]{ Recall that Ted presented the first serious attempt at resolving \textbf{P} versus \textbf{NP} \cite{wo12}.}

\footnotetext[3]{TSP polytope LP models are often modelled so that Hamiltonian graphs yield higher optimal values at tour extrema compared to non-tour extrema, as is the case here.}

\normalsize	
\section{Through the Years} Ted was well known for his work on the 4CT and many other accomplishments \cite{pfers85,gsl17,sw4ct,sw4ctFA,sw4ctSG} and as cited in \cite{swsey} and others. But you can look these things up for yourself on Google. What I want people to know is that Ted was a kind and large hearted man with a humanistic-encouraging personality. He was very creative, intelligent (an understatement), very hard working, and he had a natural sense of how to encourage people to grow and achieve their own goals. 

I first met Ted in 1983, my last year as an undergraduate in the Department of Mathematics and Statistics at the University of Guelph. Ted taught me Graph Theory. I applied to graduate school and became Ted's Masters student around 1985/1986 (more on this later). This was about the same time Ted presented his infamous (proposed) compact LP formulation of the TSP (Travelling Salesperson Problem) (\textbf{NP-hard}). Ted claimed that \textbf{P}=\textbf{NP} i.e. feasibility $\Leftrightarrow$ a graph is Hamiltonian. You can imagine the excitement and flurry of debate. Recall that deciding Hamiltonicity of a graph is \textbf{NP-complete} \cite{gjt76} and that LP had recently been proven to belong in \textbf{P} (Khachiyan 1979, and Karmarkar 1984). So Ted proposed his LP as a polynomial time algorithm to solve the TSP. For Ted's argument to hold, feasibility needed to imply the existence of integer extrema. It didn't matter that fractional extrema existed, but rather that integer extrema must also exist. A great deal of excitement and research followed.

The feasible region that Ted modelled was a polytope with both fractional and integer extrema, where the value of an objective function in a feasible solution to his LP formulation can be the same for both. Feasibility and optimality can be obtained even when no integer extrema exist. Ted continued to re-work his formulation to exclude the possibility of a feasible region whose complete set of extrema were fractional i.e. to exclude the possibility that feasibility may not imply existence of integer extrema. But in 1988, Yannakakis \cite{yanna88} proved there cannot exist any compact symmetric LP formulation of the TSP with these properties. Forever optimistic, I recall Ted commenting on the possibility of creating a compact formulation of the TSP polytope, although this was also proven impossible in 2012 \cite{fmptw12}. There's an immense rich history and cast of researchers concerning studies of the TSP,  too numerous to list here. For some early names see \cite{TSPpre}.

After graduating (M.Sc.) in 1988, I worked on a `uniquely Guelph project'. Oddly this was Ted's doing. Ted was the Director of Computing and Communication Services at that time, and he had negotiated a contract with Byte Magazine and the University of Guelph to develop a conferencing system called CoSy. But I eventually lost interest and left for work in Ottawa. Chance or destiny, the chair of the Mathematics and Statistics Department (William Smith) called me in the summer of 1990 and asked me to apply to Guelph's new Ph.D. program. So I returned to Guelph, joined Ted as his Ph.D. graduate student and never looked back.

It was a fabulous opportunity, to be working (part time) on one the world's hardest problems, and most importantly, with a professor I trusted, respected, was well published, and who committed to mentoring me. Ted encouraged me to believe in my ideas, and to go about proving an idea to be true as if it were true. He believed you could gain new insights that otherwise might not be seen. ``There's no promised outcome." he would say. ``It just helps."

\begin{figure}
\centering
\includegraphics[width=110mm,height=95mm]{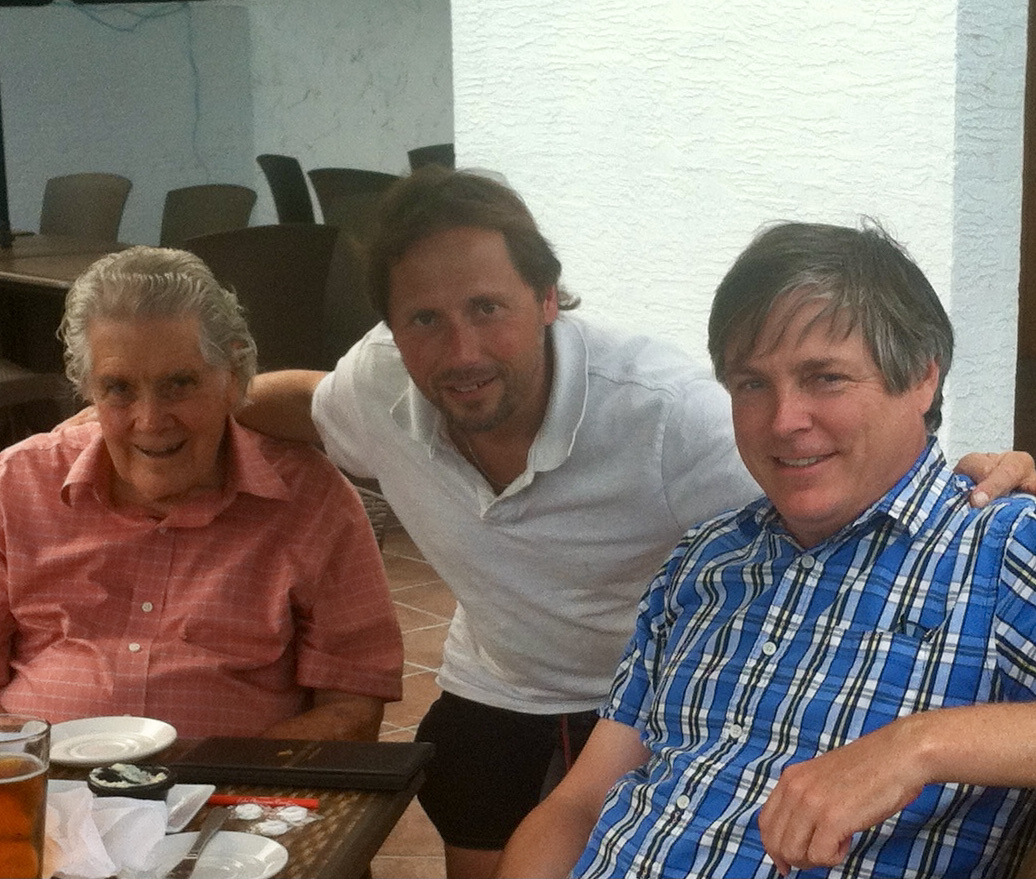}
\caption{Ted, Nic and I working together.}
\label{TedSteveNic}
\end{figure}

Our friendship started in my first graduate degree, and it blossomed during the Ph.D. We travelled together often, and attended many conferences and talks together: Baton Rouge (LSU), Boca Raton (FAU), Winnipeg and Saskatoon. Ted introduced me to many researchers e.g. I remember meeting Ralph Stanton over a plate of perogies in Winnipeg. Ted had also worked at Waterloo University in Ontario, and we'd regularly drive over (from Guelph) to the Davis Centre for more talks. We got along so well, and all the time, he was laying a  path for my own career. On the social front, Ted and Diana would celebrate my childrens' birthdays with my family, and my family would visit them at their home for Christmas parties etc. I became friends with their children as well, and we're still good friends today. In fact Ted's son Nicholas co-authored a paper with Ted, myself and others \cite{sgsb16}. I had written code to test and implement one of Ted's ideas, and Nic ran these programs on large computer farms. I remember working with Ted and Nic in 2014, at a resort in Kelowna.  We always found ways to have fun at the same time. See Figure \ref{TedSteveNic}. 

If I were to say that Ted taught me only one thing however,  it would not be mathematics. For whatever reason during my times as a graduate student, I seemed to need to prove myself. I was lucky to have met Ted because he sensed this. Ted was patient and non-confrontational in every way. Whenever I was wrong or slow to implement ideas, he would simply encourage me to believe in myself and to stop wasting time on useless insecurities. I never met any other professor like him.

Ted held strong opinions not commonly shared in his field. But he chastised no-one. I wish that all of us could say the same. He was undaunted by minor set-backs, and regularly stood out on a limb just to prove a point. He was totally secure in his own mind and never cared to waste time on positioning himself in terms of gaining approval from others - sometimes to his own detriment. I never told him about the numerous times people would ask me (conferences and talks) - in astonishment - `Does he really believe \textbf{P} = \textbf{NP}?' I only ever managed to say that no-one knows for sure. But that's not the point. I never understood the jockeying for position and Ted's approval when his ideas stood strong and unchallenged, and then the obvious distancing and abandonment Ted must have felt afterwards. I saw secret delight in the eyes of Ted's detractors when he was proven wrong. I surely didn't expect professionals to teach me these kinds of lessons. But then I learned more about Ted. I learned that he simply refused to spend time on unchangeable opinions and attitudes of others. He knew the bigger picture, that none of us have the luxury of a second life. Who can afford to cheat themselves and hold back what they want to believe and what they need to express? Gracefully accept responsibility for yourself and move on to the next task.

After Ted retired, I visited him and his family often. I can never forget Ted's 75$^{th}$ birthday. I stayed at Nic's home (the party venue), intending to help with preparations. Lots of family, friends and university colleagues were invited. But the forest fires of Kelowna (2003) were raging that same weekend. So the Swarts decided to relocate the party - ``just in case". It was wonderful i.e. the Swarts always provided for everyone, with style and class. Ted's whole family attended, coming from as far away as South Africa. Many old colleagues and friends also came to celebrate with Ted e.g. Murray Alexander, Heinz Bauschke - to mention just a few. Ted was thrilled! But rumours of evacuation, and that homes might be lost to fire overshadowed the whole event. I recall how ash and even hot embers accumulated on my rental car. When I wanted to drive, I brushed it off as if were snow. The sun looked orange (when we could see it). But usually it was just a faint glow obscured by thick smoke. Then I heard all flights in and out of Kelowna had been cancelled and I knew it was serious. Fifteen firemen almost died - trapped and surrounded by fire. We were worried for them, it was an anxious time. The whole city was on edge. Then it happened. Although Ted and Diana were already hosting a full house of out-of-town guests, 25,000 Kelowna residents were ordered to evacuate their homes with nowhere to go, including Nic and his family, and myself. I remember how Ted and Diana opened their home, and cheerfully provided for as many as possible. Nic's family and I joined others at Ted and Diana's home, and the house filled fast. So I moved to a neighbour's home and later left for a hotel, although I continued to spend my days with Ted's family, guests and evacuees. There were so many of us by then, and we shared many meals and stories with one another. Some of us bought toys for the children of families that stayed with the Swarts - the families didn't have time to pack anything. At that time, no-one would have known that two of these childeren lost all of their toys to fire. That is, we had all hoped that Kelowna would be spared and that soon everything would return to normal. But in truth it was a community tragedy in progress - 238 homes were destroyed, including the home of one of the families that stayed with us. The main point I want to share however is the kindness that Ted and his family regularly showed toward others. Respect is part of it. Stepping up to the plate, taking chances and doing the right thing even when you've reached capacity i.e. sharing your home, that's how you make it count. Ted and Diana led by example.

Looking back, Ted and his family changed my life. There are so many more stories I could share. Suffice to say, Ted's steadfast commitment and endless encouragement built me into the person I am today. He helped me in so many more ways than I can say here e.g. finance, fatherly talks and time spent together. Thank you Ted. Your efforts live on and I'll miss you.

\section{An Academic Contribution in Honour of Ted} The following sections develop definitions, a key idea, and transformation of extrema of the Birkhoff polytope leading to the definition of Ted's polytope in Section \ref{properties}. An LP formulation of the Hamilton Tour decision problem is then presented, where the feasible region is a non-compact extended formulation of Ted's polytope. A second LP formulation is presented where the feasible region is Ted's polytope embedded in the space in which it is defined. A unique external representation is computed (pure novelty) for an instance of Ted's polytope, and summary details are reported.\

\subsection{Definitions and Key Idea}\label{idea}

All graphs are simple and directed $n$-vertex graphs, interchangeably referenced by their adjacency matrices.  Definitions and the main ideas now follow.\

\begin{Defn}$P$ is the set of all $n!$ $n_xn$ 0/1 permutation matrices (also graphs / adjacency matrices), interchangably referenced as single permutation matrices whose components are $P[i,j]$, depending on context.\end{Defn}

\begin{Defn}$B_n$ =  conv$(P)$, the well known Birkhoff polytope.\end{Defn}

When permutation matrices are interpreted as adjacency matrices of $n_xn$ 2-factor graphs, $B_n$ is viewed as  the feasible region of an LP that maximizes $<P,G>$ in search of 2-factors in $G$, a problem in \textbf{P}. Its compact formulation is the set of assignment constraints.

\begin{Defn}$P_t$ is the set of all $(n-1)!$ $n_xn$ 0/1 irreflexive connected 2-factor permutations matrices (also graphs / adjacency matrices), interchangably referenced as single permutation matrices whose components are $P_t[i,j]$, depending on context. When interpreted as graphs, these are Hamilton tours, or more simply tours.\end{Defn}

\begin{Defn}$TSP_n$ =  conv$(P_t)$, the infamous TSP polytope. \end{Defn}

Polytope $TSP_n$ is the feasible region of an IP for the well known \textbf{NP-hard} Traveling Salesperson Problem, and its corresponding feasibility problem that decides existence of a tour (\textbf{NP-complete}). No-one knows exactly (although William Cook is one of the best) how to work with the TSP polytope, other than it's the solution set of an intractable set of constraints, embedded in the image space of an orthogonal projection of a convex combination of its extrema. Unlike $B_n$, $TSP_n$ has no compact formulation \cite{fmptw12}. New kinds of facets emerge for increasing $n$, and not even the number of facets are known (as a function in $n$) in general. See \cite{cook2018,cookke} for applications, experiments and many excellent studies.  

\begin{Defn}$P_{nt}$ is the set of all $n!-(n-1)!$ $n_xn$ 0/1 disconnected 2-factor permutations matrices (also graphs / adjacency matrices), the complement of $P_t$ with respect to $P$, interchangably referenced as single permutation matrices whose components are $P_{nt}[i,j]$, depending on context. When interpreted as graphs, these are simply called non-tours.\end{Defn}

\begin{rema}The union of all tours and non-tours is the set of all 2-factor graphs. i.e.  $P$=$P_t\cup P_{nt}$, and $P_t\cap P_{nt}=\emptyset$. \end{rema}

The idea that polyhedral structure is related to problem complexity is suggested by Edmonds and others as cited in \cite{cook2018}, and is best illustrated by consequences of the separation theorem \cite{GLS81,Tar86}. Thinking about how this idea might present itself with respect to the set of permutations, perhaps they provide perfect symmetry that enables the structure of $B_n$ to be expressed as the solution set of a compact formulation, while the subset of tours $P_t$ provides perfect asymmetry that enables the structure of $TSP_n$ to never possess any compact formulation. More thoughts about these two polytopes are presented in \cite{sgsb16}.\

\begin{figure}
	\centering
	\includegraphics[width=66mm,height=58mm]{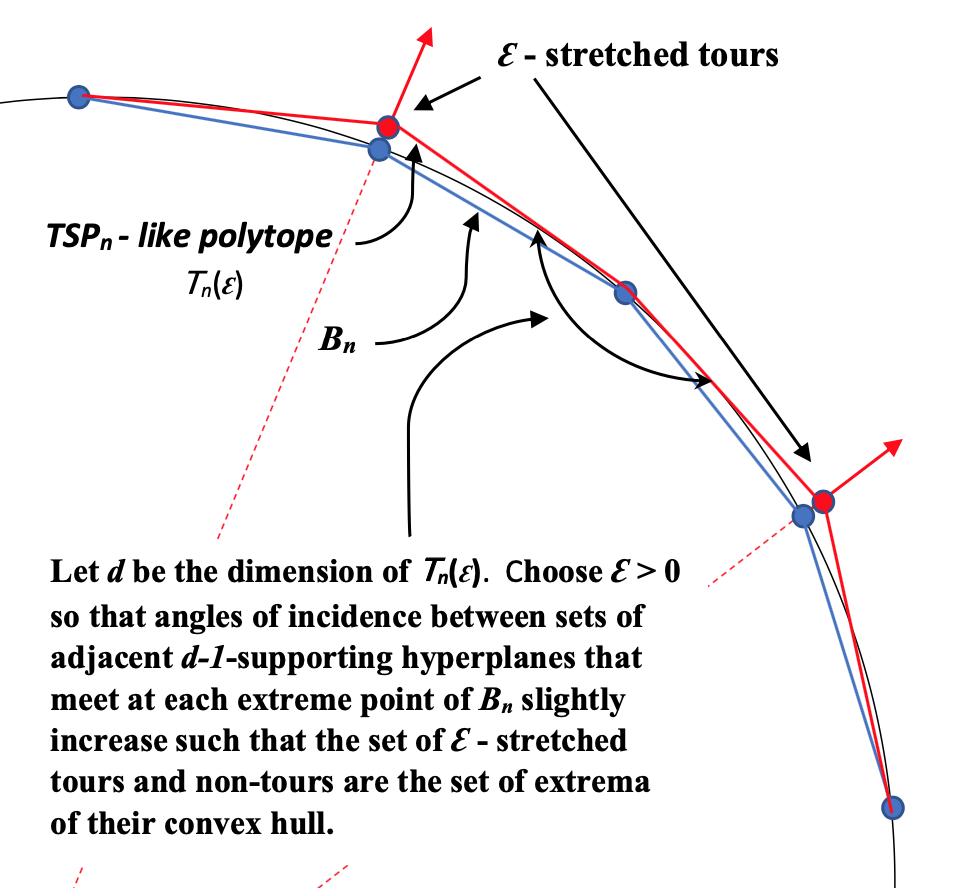}
	\caption{Key Idea:  Crafting a $TSP_n$-like Polytope via $\epsilon$-Stretched Tours of $B_n$.}
	\label{Betan}
\end{figure}

 In earlier work \cite{gismo98}, an idea was presented where, starting with a cube, extrema of its dual are crafted (from within the cube) via $\epsilon$-stretches ($\epsilon>0$) of arithmetic centers of faces of the cube away from its center. They `pop' out of faces via small $\epsilon$-stretches.  The solution set of the $2n$ constraints, a cube, AND the convex hull of the set of $2n$  extrema of its dual are then embedded in two disjoint orthogonal subspaces. The convex hull of a cube and its dual is created as a convex combination of these sets projected into the same subspace (and onto/defining the convex hull), and has a compact formulation. For small $\epsilon$ $>$ 0, extema of the dual remain extreme with respect to the convex hull of both sets. These thoughts lead to the key idea of similarly crafting an $\epsilon$-stretched $TSP_n$ from within $B_n$ by $\epsilon$-stretches of tours away from the arithmetic center of $B_n$ i.e. the set of extrema of $TSP_n$ is already a subset of the set of extrema of $B_n$, but not yet stretched. Perhaps it's possible to craft a $TSP_n$-like polytope i.e. the convex hull of 1) stretched (tour) extrema of $B_n$, and, 2) the complement set of unstretched (non-tour) extrema - and that maybe stretched tours and unstretched non-tours might remain extreme with respect to their combined convex hull. The trick, if possible, is to adjust / tune $\epsilon$ in just the right way. The goal for Ted's polytope then is to build-in a property that distinguishes extreme point tours from extreme point non-tours. An analogy is presented in Figure \ref{Betan}. Originally there are five extrema on the circle. Two $\epsilon$-stretched tours (extrema) are created via stretching two original extrema of $B_n$ away from the center of the circle. All points remain extreme, and the two $\epsilon$-stretched tours are distinguished (formalized in Lemma \ref{lem}) from the remaining set of three extrema i.e. in a circle of greater radius, although for Ted's polytope the idea is to distinguish tours from non-tours in a linear way.

\subsection{Transformation of Extrema of Polytope $\boldsymbol{B_n}$}\label{transform}
All members in $P_t$ are first transformed into \textit{stretched tours}, followed by a translation into the positive orthant. Members in $P_{nt}$ are simply translated into the positive orthant.\

\begin{Defn}\label{CentreBn}
	Define  $\overline{B_n}$ as the arithmetic center point  of $B_n$ given by\
	
\small
$\overline{B_n} =\left(\begin{array}{cccc}
	\frac{1}{n} & \frac{1}{n} & ... & \frac{1}{n}\\
	\frac{1}{n} & \frac{1}{n} & ... & \frac{1}{n}\\
	... & ... & ... & ...\\
	\frac{1}{n} & \frac{1}{n} & ... & \,\,\frac{1}{n}\
\end{array}\right)$
\normalsize
\end{Defn}
\vspace{3mm}

\textbf{Step 1. Transformation of $P_t$ Into Stretched Tours $P_t'$.} Transform all $(n-1)!$ tours, extreme points of $B_n$, to create stretched tours in the direction away from $\overline{B_n}$ i.e. ${P_t}'$= $P_t$+$\epsilon$ ($P_t$ - $\overline{B_n}$), where $\epsilon$ $>$ 0. Observe that $P_t[i,j]=0$ becomes ${P_t}'[i,j]$=$-\frac{\epsilon}{n}$, and that $P_t[i,j]=1$ becomes ${P_t}'[i,j]$=1+$\epsilon$ $-\frac{\epsilon}{n}$. See the example below, where $P_t$ is the tour $(1-2-3-4)$.\\

\begin{figure}[h!]
\begin{center}
\small
$\left(\begin{array}{cccc}
            0 & 1 & 0 & 0\\
            0 & 0 & 1 & 0\\
            0 & 0 & 0 & 1\\
            1 & 0 & 0 & 0\
\end{array}\right)\Rightarrow
\left(\begin{array}{cccc}
            0 & 1+\epsilon & 0 & 0\\
            0 & 0 & 1+\epsilon & 0\\
            0 & 0 & 0 & 1+\epsilon\\
            1+\epsilon & 0 & 0 & 0\
\end{array}\right)-\frac{\epsilon}{n}
\left(\begin{array}{cccc}
            1 & 1 & 1 & 1\\
            1 & 1 & 1 & 1\\
            1 & 1 & 1 & 1\\
            1 & 1 & 1 & 1\
\end{array}\right)$

\normalsize
\caption{How to Make A Stretched Tour.}
\label{SampeTour}
\end{center}
\end{figure}

\textbf{Step 2. Translation of Stretched Tours $P_t'$ and Non-Tours $P_{nt}$ Into the Positive Orthant.} Translate $P_t'$ into the positive orthant by an amount $\frac{\epsilon}{n}$ and rename members $Q_t$=${P_{t}}'$ +$\epsilon \overline{B_n}$=$P_t$+$\epsilon$ ($P_t$ - $\overline{B_n}$)+$\epsilon \overline{B_n}$=$P_t(1+\epsilon)$. Likewise translate and rename all $P_{nt}$ to $Q_{nt}$=$P_{nt}$ + $\epsilon \overline{B_n}$. See Figure \ref{QtQnt} where ${P_t}'$ is the stretched tour $(1-2-3-4)$, and $P_{nt}$ is the non-tour $(1-1)(2-2)(3-3)(4-4)$.

\begin{figure}[h!]
\begin{center}
\small
$\left(\begin{array}{cccc}
            0 & 1 & 0 & 0\\
            0 & 0 & 1 & 0\\
            0 & 0 & 0 & 1\\
            1 & 0 & 0 & 0\
\end{array}\right)\Rightarrow
\left(\begin{array}{cccc}
            0 & 1+\epsilon & 0 & 0\\
            0 & 0 & 1+\epsilon & 0\\
            0 & 0 & 0 & 1+\epsilon\\
            1+\epsilon & 0 & 0 & 0\
\end{array}\right)$

\vspace{4mm}

$\left(\begin{array}{cccc}
	1 & 0 & 0 & 0\\
	0 & 1 & 0 & 0\\
	0 & 0 & 1 & 0\\
	0 & 0 & 0 & 1\
\end{array}\right)\Rightarrow
\left(\begin{array}{cccc}
	1+\frac{\epsilon}{n} & \frac{\epsilon}{n} & \frac{\epsilon}{n} & \frac{\epsilon}{n}\\
	\frac{\epsilon}{n} & 1+\frac{\epsilon}{n} & \frac{\epsilon}{n} & \frac{\epsilon}{n}\\
	\frac{\epsilon}{n} & \frac{\epsilon}{n} & 1+\frac{\epsilon}{n} & \frac{\epsilon}{n}\\
	\frac{\epsilon}{n} & \frac{\epsilon}{n} & \frac{\epsilon}{n} & 1+ \frac{\epsilon}{n}\
\end{array}\right)$

\normalsize
\caption{$Q_t$=$P_t(1+\epsilon)$, and $Q_{nt}$=$P_{nt}$ + $\epsilon \overline{B_n}$.}
\label{QtQnt}
\end{center}
\end{figure}

\begin{lemm}\label{dimbetan}
dim$($aff$(Q_t,Q_{nt}))=n^2-2n+1$.
\end{lemm}\

\textbf{Proof.} Aff$(Q_t,Q_{nt})$ is a translation of aff$(B_n)$. It's well known that dim$($aff$(B_n))=n^2-2n+1$. Hence dim$($aff$(Q_t,Q_{nt}))=n^2-2n+1$. In more detail, Step 1 can be viewed as 1) translating all $P$ by $\overline{B_n}$ i.e. moving $B_n$ so its centre is the origin, 2) multiplying only the translated members of $P_t$ by a scaler multiple, followed by 3) translation back from the origin.  Since the dimension of the span of all translated extrema is unaffected by any non-zero scaling of any subset of these vectors, dim$((P_t-\frac{1}{n})(1+\epsilon),P_{nt}-\frac{1}{n})$ = dim$(P_t-\frac{1}{n},P_{nt}-\frac{1}{n})$. But dim$(P_t-\frac{1}{n},P_{nt}-\frac{1}{n})$ =  dim$(P-\frac{1}{n})$=dim$($aff$((P))$=dim$($aff$(B_n))$ = $n^2-2n+1$ i.e. dimension is preserved under translation back from the origin, and Step 2 is merely a translation ($\frac{\epsilon}{n}$) of the affine subspace created in Step 1  i.e. aff$(Q_t,Q_{nt})$. \hfill $\square$\

\subsection{Properties of Points in $\boldsymbol{Q_t}$ and $\boldsymbol{Q_{nt}}$, and Relationships to Graphs}\label{QTQ}
Let $Q$ be the set of all transformed and translated extrema of $B_n$ i.e. all points in $Q_t$ and $Q_{nt}$. Members of $Q_t$ and $Q_{nt}$ are referred to as points, and not extrema i.e.  it's not yet been shown there exists $\epsilon > 0$ such the set of points of $Q_t$ and $Q_{nt}$ is the set of extrema of conv($Q_t$,$Q_{nt}$).\

\begin{rema}Each point in $Q$ satisfies the set of scaled doubly stochastic constraints (by design).

     $\,\,\,\,\,\,\,\,\,\,\,\,\,\,\,\,\sum_{j}Q[i,j]=1+\epsilon, i=1,2,...,n$\
      
     $\,\,\,\,\,\,\,\,\,\,\,\,\,\,\,\,\sum_{i}Q[i,j]=1+ \epsilon, j=1,2,...,n$\
     
     $\,\,\,\,\,\,\,\,\,\,\,\,\,\,\,\,\,q(i,j)\in [0,1]$\
\end{rema}

\begin{rema}\label{qte}
Each point in $Q_t$ is extreme with respect to conv$(Q_t)$ i.e. each $Q_t$ lies in the surface of a strictly convex hypersphere of radius $(1+\epsilon)\sqrt{n}$, and is irreflexive.
\end{rema}      
       
\begin{rema}\label{qtne}
Each point in $Q_{nt}$ is extreme with respect to conv$(Q_{nt})$ i.e. each $Q_{nt}$ lies in the surface of a strictly convex hypersphere of radius $\sqrt{n(1+\frac{\epsilon}{n})^2+(n^2-n){\frac{\epsilon^2}{n^2}}\big)}$=$\sqrt{1+\epsilon\big(\frac{2+\epsilon}{n}\big)}$$\sqrt{n}$, and is reflexive.
\end{rema}   

\noindent\begin{theo}\label{thm1} Consider $max<G_H,Q>$ over the finite set of $|Q|=n!$ points where graph $G_H$ is Hamiltonian. Then
\[   
max<G_H,Q> is
     \begin{cases}
       \text{$=n(1+\epsilon)$, attained for a $Q$ $\in$ $Q_t$}\\
       \text{$\le n(1+\epsilon)-\epsilon$, may be attained for a $Q$ $\in$ $Q_{nt}$}\\
     \end{cases}
\]
\end{theo}

\textbf{Proof.} First consider $max<G_H,Q_t>$. For each point in $Q_t$, there are $n$ off-diagonal components at $1+\epsilon$ level, and the remaining components are at zero level. Each of these points cover a tour in correspondence with each of $(n-1)!$ tours. Thus there exists a corresponding point in $Q_t$ that covers each tour in $G_H$, and for each point, the inner product sum is therefore $n(1+\epsilon)$. Next consider $max<G_H,Q_{nt}>$. $G_H$ is simple and hence the diagonal components are at zero level. To maximize $<G_H,Q_{nt}>$, choose $G_H$ to be complete (so that there are at most $n^2-n$ off-diagonal components of $G_H$ at unit level), and consider only those $Q_{nt}$ for which diagonal components are at $\frac{\epsilon}{n}$ level, implying that for each point, there are a  maximum number of occurrences ($n$) of off-diagonal components at $1+\frac{\epsilon}{n}$ level, and $n^2-2n$ off-diagonal components at $\frac{\epsilon}{n}$ level. The maximum possible value of $<G_H,Q_{nt}>$=$n(1+\frac{\epsilon}{n})$ + $(n^2-2n)(\frac{\epsilon}{n})$=$n(1+\epsilon)-\epsilon$.\hfill $\square$\\

\noindent\begin{theo}\label{thm2} Consider $max<G_{NH},Q>$ over the finite set of $|Q|=n!$ points where graph $G_{NH}$ is non-Hamiltonian graph. Then
\[   
max<G_{NH},Q>= 
     \begin{cases}
       \text{$\le (n-1)(1+\epsilon)$, may be attained for a $Q$ $\in$ $Q_t$}\\
       \text{$\ll n(1+\epsilon)-\epsilon$, for all $Q$ $\in$ $Q_{nt}$}\\
     \end{cases}
\]
\end{theo}

\textbf{Proof.} First consider $max<G_{NH},Q_t>$. Recall from above that for each point in $Q_t$, there are $n$ off-diagonal components at $1+\epsilon$ level, and the remaining components are at zero level. As before, each of these points cover a tour in correspondence with each of $(n-1)!$ tours. Since $G_{NH}$ has no tour, no point in $Q_t$ covers a tour in $G_{NH}$, and the inner product sum is at most $(n-1)(1+\epsilon)$. Next, consider $max<G_{NH},Q_{nt}>$. Since $G_{NH}$ is simple, the diagonal components are again at zero level. To maximize $<G_{NH},Q_{nt}>$, choose $G_{NH}$ to be critically non-Hamiltonian with as many possible arcs at unit level i.e. unlike above, $G_{NH}$ cannot be complete otherwise it becomes Hamiltonian. Suppose there are $f(n)<n^2-n$ components of $G_{NH}$ at unit level. Now consider only those $Q_{nt}$ for which diagonal components are at $\frac{\epsilon}{n}$ level so that for each point there are a maximum number of occurrences ($n$) of off-diagonal components at $1+\frac{\epsilon}{n}$ level, and $n^2-2n$ off-diagonal components at $\frac{\epsilon}{n}$ level. The maximum possible value of $<G_{NH},Q_{nt}>$ = $n(1+\frac{\epsilon}{n})$ + $f(n)(\frac{\epsilon}{n})$ $\ll$ $n(1+\epsilon)-\epsilon$.\hfill $\square$\

\begin{rema}
Observe that $Q_t$ are irreflexive and that $Q_{nt}$ are reflexive, leading to the existence of many hyperplanes that separate conv($Q_t$) and conv($Q_{nt}$). This separation property is the foundation reason behind proofs of Theorems \ref{thm1} and \ref{thm2}.
\end{rema}

\begin{lemm}\label{lem}Consider $max<G,Q>$ over the finite set of $|Q|=n!$ points. Then $max<G,Q>$=$n(1+\epsilon)$ $\Leftrightarrow$ graph $G$ is Hamiltonian.
\end{lemm}

\subsection{Ted's Polytope}\label{properties}\

\begin{Defn} \label{beta}
$\mathcal{T}_n(\epsilon)$ = conv$(Q_t,Q_{nt})$ where $\epsilon$ $>$ 0, Ted's polytope.
\end{Defn}

Recall Lemma \ref{dimbetan}, dim$($aff$(Q_t,Q_{nt}))$=dim$($aff$(\mathcal{T}_n(\epsilon)))=n^2-2n+1$, and following from Remark \ref{qte}, $\mathcal{T}_n(\epsilon)$ has at least $|Q_t|$ extrema. 

\subsubsection{Extrema of $\boldsymbol{\mathcal{T}_n(\epsilon)}$, $\boldsymbol{\epsilon_g}$ and $\boldsymbol{\epsilon_{max}}$}\label{extrema}\

\begin{Defn}Given $n$, $\epsilon$, and $\mathcal{T}_n(\epsilon)$, suppose all members of $Q$ are precisely the set of extrema of $\mathcal{T}_n(\epsilon)$. Then  $\epsilon$ is said to be good, called $\epsilon_g$.
\end{Defn}

\begin{Defn}Given $n$, $\epsilon_g$, and $\mathcal{T}_n(\epsilon_g)$, the greatest $\epsilon_g$ is $\epsilon_{max}$.
\end{Defn}

An LP model was created, coded and run so that for $n$=4, 5 and 6, it determined the set of extrema of $\mathcal{T}_n(\epsilon)$ as a function of chosen $\epsilon$. For each $Q_i$ in $Q_{nt}$, noting that the set of $Q_t$ are extreme, the LP attempted  to create $Q_i$ as a convex combination of members in $Q_t$ union $Q_{nt}$, not including $Q_i$. If every $Q_i$ is not a convex combination of these members, then all members in $Q_t$ and $Q_{nt}$ are extreme with respect to  $\mathcal{T}_n(\epsilon)$. \

Results are shown in Table \ref{HCBETANX} for $\epsilon$=1, 5, 10 and 20, yielding existence of an $\epsilon_g$ for each of $n$=4, 5 and 6. Entries in the table are counts\footnotemark[8] of extrema of $\mathcal{T}_n(\epsilon)$ written in the form 1 - 2 - 3, and classified as;
\begin{enumerate}
\item{transformed, irreflexive tours i.e. irreflexive $Q_t$, stretched tours}
\item{translated, irreflexive non-tours i.e. reflexive $Q_{nt}$}
\item{translated, reflexive non-tours i.e. reflexive $Q_{nt}$}
\end{enumerate}
Unbolded entries indicate that the set of $Q_t$ union $Q_{nt}$ is the set of extrema of $\mathcal{T}_n(\epsilon)$. As $\epsilon$ increases, the convex hull of irreflexive $Q_t$ corresponding to $\epsilon$-stretched tours expands, and so too does all of $\mathcal{T}_n(\epsilon)$ as all $Q_{nt}$ are further translated away from the origin and further into the positive orthant. Consequently some members of $Q_{nt}$ become convex combinations of irreflexive $Q_t$ and other members of $Q_{nt}$, and are no longer extrema of $\mathcal{T}_n(\epsilon)$. Interestingly and easy to show, only those $Q_{nt}$ corresponding to translated irreflexive non-tours have the potential to lie interior to $\mathcal{T}_n(\epsilon)$, and, bolded entries in the table indicate when this happens.\

For $n$=4, no $Q_{nt}$ corresponding to irreflexive non-tours fall interior to $\mathcal{T}_4(\epsilon)$ (for all choices of $\epsilon$). For $n$=5, all $Q_{nt}$ corresponding to irreflexive non-tours fall interior to $\mathcal{T}_5(\epsilon)$ for $\epsilon \ge 5$. For $n$=6, some $Q_{nt}$ corresponding to irreflexive non-tours fall interior to $\mathcal{T}_6(\epsilon)$ starting when $\epsilon$ = 5, and more  $Q_{nt}$ continue to fall interior to $\mathcal{T}_6(\epsilon)$ as $\epsilon$ increases beyond 5.\

\footnotetext[2]{ Maybe there exists sets $E$ of individually parameterized $\epsilon$, specific to individual tour stretches, that together create an ideal symmetry in $\mathcal{T}_n(E)$. That is, as the parameter is adjusted past the point where $\mathcal{T}_n(E)$ achieves its ideal symmetry, an instant and complete transition might then occur. Lemma \ref{lem} might change to read something like ``...  $max<G,Q>$ $\ge$ $n(1+E_{min})$ $\Leftrightarrow$ graph $G$ ..."}

\begin{table}[h!] 
\centering
\caption{Classification and Counts of Extrema of $\mathcal{T}_n(\epsilon)$ as a Function of $n$ and $\epsilon$.}\label{HCBETANX} 
	{ \footnotesize
		\begin{tabular}{c|c|c|c|c|c c}  
			                      &                                  &                                                                         &                                                                          &                                                                  &                             \\
    		$n$               &   $B_n$	                & $\mathcal{T}_n(1)$                      & $\mathcal{T}_n(5)$                     & $\mathcal{T}_n(10)$                & $\mathcal{T}_n(20)$ \\    
    		    			 \hline  
    		            
		    4                    &   6-3-15                 &    6-3-15                                                     &    6-3-15                                                         & 6-3-15                                                     &  6-3-15            \\
			5                    &    24-20-76           &   24-20-76                                               &    \textbf{24-0-76}                                      &  \textbf{24-0-76}                                  &  \textbf{24-0-76}  \\
			6                     &    120-145-455   &   120-145-455                                        &  \textbf{120-26-455}                                   &  \textbf{120-19-455}                          &  \textbf{120-13-455}   \\
					[1ex]    
		\end{tabular} 
    }
\end{table}

These results are fascinating on there own. No matter, $\mathcal{T}_n(\epsilon_g)$ exists for $n$=4, 5, and 6 ($\epsilon_g=1$ is sufficient), leading to Conjecture \ref{con1}.

\footnotetext[8]{These counts varried depending on machine zero between 10$^{-8}$ through 10$^{-10}$.}

\begin{Prop}\label{con1}
	Polytope $\mathcal{T}_n(\epsilon_g)$ exists for all $n \ge 4$.
\end{Prop}

\subsubsection{An Extreme Point Formulation of $\boldsymbol{\mathcal{T}_n(\epsilon_g)}$}\label{external} 
Given an $\epsilon_g$, write each point in $Q_t$ as an $n^2$ long unbolded column vector $q_t$, whose components are $Q_t[i,j]$ as read from left to right and top to bottom. There are $(n-1)!$ \textbf{$q_t$}. Similarly create the set of $n!-(n-1)!$ $q_{nt}$. Define $\boldsymbol{\alpha}$ as a non-negative column vector of variable scalars that sum to 1, and define $n^2$ long column vector variable \textit{\textbf{q}} to be convex combinations of members from $q_t$ and $q_{nt}$ depending upon how $\boldsymbol{\alpha}$ is assigned. See Figure \ref{extrempt}.\

An extreme point formulation of $\mathcal{T}_n(\epsilon_g)$ is an extended formulation, and has size $\mathcal{O}(n!)$, and its set of extrema are embedded in $\mathcal{R}^{n^2}$ coded by \textbf{\textit{q}}. For comparison, $B_n$ is modelled by the assignment constraints and has a compact formulation of size $\mathcal{O}(n^2)$. Its set of extrema are also embedded in $\mathcal{R}^{n^2}$. Both polytopes possess unique representations in $\mathcal{R}^{n^2-2n+1}$ i.e. for $B_n$ there are $n^2$ unique intersecting half-spaces, defined by $n^2$ unique facet-inducing inequalities. A natural question now arises: `How many unique intersecting half-spaces define $\mathcal{T}_n(\epsilon_g)$?' A small discussion about unique representations is presented in Section \ref{image}, including an example. 

\begin{figure}[h!]
	\begin{center}
		\small
		$\left(\begin{array}{c|c}
			q_t\,\,\,\,\,    q_{nt} & -I\
			
			\\\hline    
			\textbf{1}  & \textbf{0} \
		\end{array}\right)
		\left(\begin{array}{c}
			{\boldsymbol{\alpha}}\
			\\    
			\textit{\textbf{q}}\
			
		\end{array}\right)=
		\left(\begin{array}{c}
			\textbf{0}\
			\\    
			1 \
		\end{array}\right)$\\
		\vspace{2mm}
		$\boldsymbol{\alpha} \ge 0$

		\normalsize
		\caption{Extreme Point Formulation Where $\mathcal{T}_n(\epsilon_g)$ is Embedded in $R^{n^2}$.}
		\label{extrempt}
	\end{center}
\end{figure}

A non-compact LP formulation of the Hamilton Tour decision problem is presented in Figure \ref{hami}, followed by a second formulation presented in Figure \ref{hami2} where the feasible region is $\mathcal{T}_n(\epsilon_g)$ embedded in $R^{n^2}$. Let $G$ be an arbitrary graph written as an $n^2$ long column vector \textit{\textbf{g}} whose components are $G[i,j]$ as read from left to right and top to bottom.\\

\begin{figure}[h!]
	\begin{center}
		\small
		$max<$\textbf{\textit{g$^T$}},\textbf{\textit{q}}$>$ subject to;\
			\vspace{2mm}
		
		$\left(\begin{array}{c|c}
			q_t\,\,\,\,\,    q_{nt} & -I\
			
			\\\hline    
			\textbf{1}  & \textbf{0} \
		\end{array}\right)
		\left(\begin{array}{c}
			{\boldsymbol{\alpha}}\
			\\    
			\textit{\textbf{q}}\
			
		\end{array}\right)=
		\left(\begin{array}{c}
			\textbf{0}\
			\\    
			1 \
		\end{array}\right)$\\
		\vspace{2mm}
		$\boldsymbol{\alpha} \ge 0$		\
			\vspace{2mm}
		
		where $G$ is Hamiltonian if and only if $max<$\textbf{\textit{g$^T$}},\textbf{\textit{q}}$>$=$n(1+\epsilon)$.
			\vspace{2mm}
		
		\normalsize
		\caption{LP Formulation of the Hamilton Tour Decision Problem.}\label{hami}
	\end{center}
\end{figure}

\begin{figure}[h!]
	\begin{center}
		\small
		$max<$\textbf{\textit{g$^T$}},\textbf{\textit{q}}$>$ subject to;\
						\vspace{2mm}
						
		$\boldsymbol{q}$ $\in$ $\boldsymbol{\mathcal{T}_n(\epsilon_g)}$
				\vspace{2mm}
				
		$\boldsymbol{q} \ge 0$		\

			where $G$ is Hamiltonian if and only if $max<$\textbf{\textit{g$^T$}},\textbf{\textit{q}}$>$=$n(1+\epsilon)$.
					\vspace{2mm}
					
		\normalsize
		\caption{LP Formulation of the Hamilton Tour Decision Problem.}\label{hami2}
	\end{center}
\end{figure}

\subsubsection{A Unique External Representation of $\boldsymbol{\mathcal{T}_n(\epsilon_g)}$}\label{image}

In general, a (full dimensional) facet of a polytope is defined as the intersection set of a polytope and the solution set of a bounding hyperplane (facet-inducing inequality) whose affine hull is dimension dim(polytope) - 1.  See \cite{bron83} for details. An orthogonal projection of the extreme point formulation of a polytope onto an image embedded in a dim(polytope) dimensional vector space results in a minimally sized, finite unique set of irredundant facet-inducing inequalities, and these inequalities can be computed using a variety of techniques \cite{bartels,dg08,gismo98,gismo99,gismo03,gismo01,heller55}. The set of inequalities is (called) a unique external representation of the polytope.\

\textbf{Example. $\boldsymbol{\mathcal{T}_n(\epsilon_g)}$}.
A suitable set of image space variables for general $n$ is shown in Figure \ref{liindep2}. For $n$=4 and $\epsilon_g$=1, and dropping the last row and column for each of the 24 extrema (6 $Q_t$ and 18 $Q_{nt}$) of  $\mathcal{T}_4(1_g)$, a unique set of 508\footnotemark[1] irredundant facet-inducing inequalities (unique external representation) is computed whose corresponding 508 half-space interestion set of dimension 9  is $\mathcal{T}_4(1_g) \in \mathcal{R}^9$.

	\begin{figure}[h!]
	\begin{center}
			
			\small
			    $\left(\begin{array}{ccccc}
					Q[1,1]      & Q[1,2]      & ...   &      Q[1,n-1]       & X\\
					Q[2,1]      & Q[2,2]      & ...   &      Q[2,n-1]       & X\\
					...         & ...         & ...   &      ...            & X\\
					Q[n-1,1]    & Q[n-1,2]    & ...   &      Q[n-1,n-1]     & X\\
					X           &        X    & X     &      X              &\,\,\, X\
				\end{array}\right)$
			
			\caption{Image Space Variables for $\mathcal{T}_n(\epsilon_g)$ in $\mathcal{R}^{n^2-2n+1}$.} \label{liindep2}
			\normalsize
		         \end{center}
		\end{figure}



		
	

\footnotetext[1]{Computed using the Ray algorithm \cite{gismo99}. This technique is probabilistic, and can only provide a lower bound on facet counts.}

\section{Discussion} Two comments are now noted, that result from having lifted the set of tours from $B_n$, in the form of $Q_t$ as separated from $Q_{nt}$.\

By design, all $Q_{t}$ corresponding to all \textit{irreflexive tours} $P_t$ lie in the hyperplane $\sum Q_t[i,i]$ = 0. The remaining set of irrefexive non-tours are translated to a subset of $Q_{nt}$ whose $Q_{nt}[i,i]$=$\frac{\epsilon}{n}$. Hence they are likewise separated from tours and reflexive non-tours, and lie in the hyperplane $\sum Q_{nt}[i,i]$ = $\epsilon$. This is the first observation / comment.\

The second comment is more speculative / questioning. Recall from Table \ref{HCBETANX}, and how the subset of $Q_{nt}$ corresponding to \textit{irreflexive non-tours} appear to fall interior to $\mathcal{T}_n(\epsilon>\epsilon_{max})$, for increasing $\epsilon$ (the middle number decreases for increasing $\epsilon$ for $n$=6). That is, as a function of $n$ and $\epsilon$, there exist convex combinations of $Q_{nt}$ and at least one member of $Q_t$ that lie strictly interior to $\mathcal{T}_n(\epsilon>\epsilon_{max})$. Recall Remark \ref{qtne}. It's not clear what may be happening here, whether or not this is generally true, and whether or not understanding this strange consequence may or maynot be useful.\

My tribute to Ted is now complete. If he were here today, I think he'd chuckle at my `rationalization' of giving him a gift polytope. I can even imagine hearing him comment ``You're giving me a polytope ..?? Ha!". You see however, it truly is a perfect gift. This is how we met, how we started together back in the 1980's when Ted presented his original LP model, and, this is how we finished. Therefore in loving memory I dedicate $\mathcal{T}_n(\epsilon_g)$ to Ted Swart, one of my dearest, gentle and most creative friends, my colleague and my mentor.



\bibliography{ReferenceDB1}
\bibliographystyle{siam}
\end{document}